\documentclass[10pt,a4paper]{article}

\usepackage[latin1]{inputenc}
\usepackage{amsmath, amssymb}

\usepackage{amsthm}

\newtheorem{thm}{Theorem}[section]
\newtheorem{prop}[thm]{Proposition}
\newtheorem{lemma}[thm]{Lemma}

\newtheorem{conjecture}{Conjecture}

\DeclareMathOperator{\Ber}{Ber}
\DeclareMathOperator{\var}{Var}
\DeclareMathOperator{\fix}{fix}
\DeclareMathOperator{\abs}{abs}
\DeclareMathOperator{\R}{\mathbb{R}}
\DeclareMathOperator{\E}{\mathbb{E}}
\DeclareMathOperator{\Prb}{\mathbb{P}}

\begin{document}

\author{Daniel Lanoue}
\date{\today}
\title{The iPod Model}

\maketitle

\begin{abstract}
We introduce a Voter Model variant, inspired by social evolution of musical preferences. In our model, agents have preferences over a set of songs and upon meeting update their own preferences incrementally towards those of the other agents they meet. Using the spectral gap of an associated Markov chain, we give a geometry dependent result on the asymptotic consensus time of the model.
\end{abstract}

\section{Introduction}

The terminology of Finite Markov Information Exchange (FMIE) models has been introduced \cite{aldous2013} \cite{aldous2012} as a catch-all for the interpretation of Interacting Particle Systems (IPS) models as stochastic social dynamics. Many important and classical models fit under this two-level framework; the bottom level a meeting model among agents, and the top level an information exchange algorithm performed at each meeting. 

For classic IPS models, such as the Voter Model, with a simple meeting algorithm the FMIE perspective is perhaps unnecessary. In this paper however, we will introduce and study a (much) generalized Voter Model - inspired by the evolution of musical preferences among a group of friends - as an FMIE process.

\subsection{The iPod Model}

Here we introduce the iPod FMIE model. The underlying framework of the stochastic process is a weighted graphs $\mathfrak{G}$ on $N$ vertices. We will refer to each vertex as an agent and occasionally to our vertex set as $I$.  Associated to the edges are symmetric meeting rates $\nu_{i,j}$ for $1 \leq i \neq j \leq N$. We assume that all meeting rates are normalized, i.e.
$$
\sum_j \nu_{i, j} = 1
$$
for all agents $i$.

Each agent $i$ is equipped at each time $t$ with a probability measure $X_t(i)$ on $\{1, 2, \ldots, \sigma \}$ which we will reference by its distribution $X_t^{k}(i)$ for $1 \leq k \leq \sigma$.

We consider $\sigma$ as a fixed number of songs and $X_t^{k}(i)$ the preference of agent $i$ at time $t$ for song $k$. The stochastic process $X_t$ updates over time as follows. Between every pair of agents $i, j$ we associate a Poisson process with rate $\nu_{i,j}$ whose times we refer to as meetings between $i$ and $j$. At a meeting time $t$ between agents $i$ and $j$, each agent picks a song $\sigma_i$ and $\sigma_j$ independently and distributed according to $X_{t-}(i)$ and $X_{t-}(j)$.

We interpret this as each agent choosing a song to play to the other agent based on their preferences. After agent $i$ hears the song $j$ chose he updates his preferences according to
$$
X_t^{\sigma_j}(i) = (1 - \eta)X_{t-}^{\sigma_j}(i) + \eta
$$
and
$$
X_t^{k}(i) = (1 - \eta)X_{t-}^{k}(i)
$$
for all other $k \neq \sigma_j$. Here $0 < \eta < 1$ is a fixed interaction parameter. Agent $j$ updates her preferences similarly. It is immediate that if $X_{t-}(i)$ is a probability measure than so is $X_t(i)$. Note that we are implicitly working with cadlag paths.

Analogous to results on the consensus time of the Voter Model - for instance \cite{cox1989coalescing} or more generally \cite{oliveira2013} - in this paper we will estimate the fixation time (to be defined) of the iPod process. Interestingly, again similar to the Voter Model our proof will explore a connection between this process and the Wright-Fisher diffusion \cite{cox1989coalescing}.

A special feature of the model (\protect \MakeUppercase {P}roposition\nobreakspace \ref {M Martingale}) is that the average (over agents) preference for a given song evolves as a martingale, analogous to the total proportion of agents with a given opinion on the voter model. This distinguishes the iPod model from many other variants of the voter model that have been studied \cite{castellano2009}.

\subsection{Fixation Time}

We will be focused on estimating the fixation time $T_{\fix}$ of the iPod process. Every time two agents meet at least one distinct song is played between them and so at least one of the $\sigma$ songs is played infinitely often. Given that only one song is played infinitely often, we define $T_{\fix}$ to be the last time any other song is played.

We note that $T_{\fix}$ is not a stopping time and a priori could be infinite, i.e. if more than one song is played infinitely often. However, we will show that this is not the case and in fact $T_{\fix}$ has finite expectation, the bounding of which will be our primary goal.

\begin{thm}{There exists a constant $C(\eta)$ so that from any initial configuration of $\sigma$ songs, the fixation time $T_{\fix}$ has expectation
$$
\E T_{\fix} \leq C(\eta) \ln(\sigma) \frac{N}{\lambda} ,
$$
where $\lambda$ is the spectral gap of $\mathfrak{G}$.
}
\label{Main Theorem}
\end{thm}

The spectral gap $\lambda$ of reversible Markov chain is interpreted as its asymptotic rate of convergence to its stationary distribution, and can be defined by the second eigenvalue of the chain's transition matrix \cite{levin2009markov}. In our setting, we define the spectral gap $\lambda$ in terms of the edge weights $\nu_{i,j}$. First, for any function $f \colon I \rightarrow \R$ we define the Dirichlet form $\varepsilon(f, f)$ by
$$
\varepsilon(f, f) = \sum_{i, j} \frac{\nu_{i,j}}{2 N} \left(f(i) - f(j) \right)^2.
$$
The spectral gap $\lambda$ is then defined in our context by the extremal characterization
$$
\lambda = \inf_{f \colon I \rightarrow \R \vert \var(f) \neq 0} \frac{ \varepsilon(f, f)}{\var(f)}.
$$

There is extensive literature \cite{levin2009markov} giving order of magnitude bounds on the $N \rightarrow \infty$ asymptotic behaviour of $\lambda_N$ for particular families of $N$-vertex graphs.  For such families, \protect \MakeUppercase {T}heorem\nobreakspace \ref {Main Theorem} gives an order of magnitude upper bound on the asymptotic fixation time, for fixed $\eta$ and $\sigma$. We will show (\protect \MakeUppercase {T}heorem\nobreakspace \ref {Complete Graph Lower Bound}) the tightness of this bound in the case of particular special family of graphs.

\section{Projection on a Single Song}

Our main technique will be focusing on the projection of our system to a single song. For some fixed (arbitrary) song $k$ we will consider only $X^k(i)$ which we will write simply as $x(i)$ dropping the $k$. When two agents $i,j$ meet, each independently chooses to either play song $k$ or not; with probability $x(i)$ and $x(j)$ respectively. Writing $\Ber(x(i))$ and $\Ber(x(j))$ for independent Bernoulli variables with given success parameters, we see that if $i$ and $j$ meet at time $t$ that
$$
 x_t(i) = (1 - \eta)x_{t-}(i) + \eta \Ber(x_{t-}(j)),
$$
with $x(j)$ updating similarly. At such a meeting, for all other agents $k \neq i, j$, $x(k)$ remains unchanged. 

This implies that the evolution of any given song can be considered separately from the others - though not independently. We will therefore focus first on the FMIE system $\{ x_t(i) \}_{i \in I, t \geq 0}$ evolving as above and then later return to the original multi-song model. The primary object of study in our one song model will be the average preference for the song, written
$$
M_t = \sum_{i \in I} \frac{x_t(i)}{N}.
$$

Our goal in this section will be bounding how long it takes $M_t$ to approach the boundary $\{0, 1 \}$. Specifically we will prove a bound on the stopping time
$$
S = \inf \{ t \geq 0 \colon M_t \notin (\frac{\eta}{2N}, 1 - \frac{\eta}{2N}) \}.
$$
We will use the shorthand $x_t = \{x_t(i) \colon 1 \leq i \leq N \}$ for the configuration at time $t$. In particular, we will often use $x_0$ for an arbitrary initial configuration. By comparison, we will use $X_t$ (respectively $X_0$) for a configuration of the multi-song model.

To state our bound, we introduce the function $\phi(x)$ given by
\begin{align}
\phi(x) = - x \ln(x) - (1 - x) \ln(1 - x).
\end{align}

\begin{thm}{There exists a constant $A(\eta)$ so that from any initial configuration $x_0$
$$
\E_{x_0} S \leq A(\eta) \frac{N}{\lambda} \phi(M_0),
$$
where $\lambda$ is the spectral gap of $\mathfrak{G}$.
}
\label{Main Escape Theorem}
\end{thm}

 To prove our theorem, we first estimate how long it takes $M_t$ to exit small intervals. Then, we use embedding to compare $M_t$ to the Wright-Fisher Diffusion.

\subsection{Derived Quantities}
We will begin by analysing a few quantities derived from $x_t$. For ease of notation we will occasionally drop $t$. Our primary object of study will be the ($L^1$) average of the preferences $x(i)$, denoted $M_t$ which is introduced above.
We will repeatedly make use of the following lemma on the step sizes of $M_t$.

\begin{lemma}{If $t$ is a meeting time then
$$
| M_t - M_{t-}| \leq \frac{2 \eta}{N}.
$$}
\label{M Steps}
\begin{proof}
If agent $i$ is involved in a meeting at $t$, then either
$$
x_t(i) = (1 - \eta)x_{t-}(i) \text{  or  } x_t(i) = (1 - \eta)x_{t-}(i) + \eta,
$$
and so
$$
| x_t(i) - x_{t-}(i) | \leq \eta.
$$
As only two agents are involved in any meeting, our bound follows easily.
\end{proof}
\end{lemma}

As a warm-up for the more complicated quantities to appear later, we begin by showing that $M_t$ evolves as a continuous time martingale. We here implicitly use the filtration $\mathfrak{F}_t$ generated by $\{x_t(i) \}_{i \in I, t \geq 0}$. Also, note that we may clearly assume that almost surely meeting times between agents are unique and that the set of meeting times has no accumulation point.

We will make use of the process dynamics notation
$$
\E \left( d A_t \vert \mathfrak{F}_{t-} \right) = (\text{resp. } \geq , \leq) B_t dt
$$
to mean that
$$
A_t - A_0 -  \int_0^t B_r dr
$$
is a martingale (respectively submartingale, supermartingale). Clearly this notation is compatible with arithmetic operations. To calculate a process's dynamics, we make repeated use of the following lemma, the proof of which is straightforward.

\begin{lemma}{Let $A_t$ be a function of the $x_t(i)$. Then
$$
\E \left( d A_t \vert \mathfrak{F}_{t-} \right) = \sum_{i,j} \nu_{i,j} \E \left(A_t - A_{t-} | i \text{ and } j \text{ meet at } t\right)dt
$$} \label{Process Dynamics}
\end{lemma}

In particular, for the average preference $M_t$ we have the following dynamics.

\begin{prop}{With respect to the filtration $\mathfrak{F}_t$, $M_t$ is a continuous time martingale.}
\label{M Martingale}
\begin{proof}
To begin we note that since $\E \Ber(x_t(j)) = x_t(j)$ we have that
$$
\E \left(x_t(i) \vert i \text{ and } j \text{ meet at time } t, \mathfrak{F}_{t-} \right) = (1 - \eta)x_{t-}(i) + \eta x_{t-}(j),
$$
and similarly for $x_t(j)$. Summing both we find that
$$
\E \left(x_t(i) + x_t(j) \vert i \text{ and } j \text{ meet at time } t, \mathfrak{F}_{t-} \right) = x_{t-}(i) + x_{t-}(j).
$$
As only $x(i)$ and $x(j)$ change at such a time $t$, this gives us that
$$
\E\left(M_t \vert i \text{ and } j \text{ meet at time } t , \mathfrak{F}_{t-} \right) = M_{t-},
$$
which clearly implies that
$$
\E \left( d M_t \vert \mathfrak{F}_{t-} \right) = 0,
$$
i.e. $M_t$ is a martingale.
\end{proof}
\end{prop}

We next look at the process dynamics of $M^2_t$. To do so we introduce the  quantity $Q_t$ given by
$$
Q_t = \sum_{i \in I} \frac{x_t(i)(1 - x_t(i))}{N}.
$$
In particular we use Lemma\nobreakspace \ref {Process Dynamics} to calculate the following.

\begin{prop}{The variation $M^2_t$ satisfies
$$
\E \left( d M^2_t \vert \mathfrak{F}_{t-} \right) = \frac{2 \eta^2}{N} Q_t dt.
$$}
\label{M2 Dynamics}
\begin{proof}
 As before, we begin by calculating  that for $k \neq i, j$, since $x(k)$ does not change after a meeting between $i$ and $j$ that:
$$
\E \left( x_t(k)(x_t(i) + x_t(j)) \vert i \text{ and } j \text{ meet at time } t , \mathfrak{F}_{t-} \right) = x_{t-}(k)(x_{t-}(i) + x_{t-}(j)).
$$
Next we calculate that
\begin{align*}
\E ( x_t^2(i) &\vert i \text{ and } j \text{ meet at } t , \mathfrak{F}_{t-} ) \\
&= (1 - \eta)^2 x_{t-}^2(i) + 2 \eta (1 - \eta) x_{t-}(i) x_{t-}(j) + \eta^2 x_{t-}(j),
\end{align*}
and similarly for $x^2(j)$. Finally we have that
\begin{align*}
\E ( x_t(i)& x_t(j) \vert i \text{ and } j \text{ meet at } t , \mathfrak{F}_{t-} ) \\
 &= (1 - \eta)^2 x_{t-}(i) x_{t-}(j) + \eta(1 - \eta)[x_{t-}^2(i) + x_{t-}^2(j)] + \eta^2 x_{t-}(i) x_{t-}(j).
\end{align*}
Putting this all together we find that
\begin{align*}
\E ( (\sum_i x_t(i))^2 &\vert i \text{ and } j \text{ meet at } t , \mathfrak{F}_{t-} ) \\
&= (\sum_i x_{t-}(i))^2 + \eta^2(x_{t-}(i) - x_{t-}^2(i) + x_{t-}(j) - x_{t-}^2(j)).
\end{align*}
Using Lemma\nobreakspace \ref {Process Dynamics}, summing over $i, j$ and normalizing by $N^2$ we find that
$$
\E \left( d M^2_t \vert \mathfrak{F}_{t-} \right) = \frac{2 \eta^2}{N} Q_t dt.
$$
\end{proof}
\end{prop}

Instead of $M^2$, we will often be more concerned with $M_t(1 - M_t)$. As $M_t$ is a martingale, from \protect \MakeUppercase {P}roposition\nobreakspace \ref {M2 Dynamics} we easily have that
$$
\E \left( d M_t(1 -M_t) \vert \mathfrak{F}_{t-} \right) = -  \frac{2 \eta^2}{N} Q_t dt.
$$ 

A central tool for the study of the underlying Markov Chain on $\mathfrak{G}$ is the Dirichlet form $\varepsilon$. We recall that the Dirchilet form $\varepsilon(f, f)$ for a function $f \colon I \rightarrow \R$ is defined as
$$
\varepsilon(f,f) = \sum_{i, j} \frac{\nu_{ij}}{2N}(f(i) - f(j))^2.
$$
We will write $\varepsilon(x_t, x_t)$ for the Dirichlet form of the function $i \mapsto x_t(i)$. 

The main fact that we will need about the Dirichlet form is its relationship to the spectral gap. We recall the definition of the spectral gap of a Markov Chain is given by
$$
\lambda = \inf_{f \colon I \rightarrow \R \vert \var(f) \neq 0} \frac{ \varepsilon(f, f)}{\var(f)},
$$
where $\var(f)$ is the variance of the function $f(i)$ with respect to the uniform measure on $I$. A simple but important fact we make repeated use of is that $0 < \lambda \leq 1$.

Following Lemma\nobreakspace \ref {Process Dynamics} we can calculate $dQ$.

\begin{prop}{The sum $Q_t$ satisfies
$$
\E \left( d Q_t | \mathfrak{F}_{t-} \right) = 4 \eta(1 - \eta) \varepsilon(x_t, x_t) dt - 2 \eta^2 Q_t dt,
$$
as well as
$$
\E \left( d Q_t | \mathfrak{F}_t \right) \geq 4 \lambda \eta(1 - \eta) M_t(1 - M_t) dt - \left( 2 \eta^2 + 4 \lambda \eta(1 - \eta) \right) Q_t dt.
$$
}
\label{Q Dynamics}
\begin{proof}
We begin by noting that $Q_t = M_t - \sum_i \frac{x_t^2(i)}{N}$ and so
$$
\E \left( d Q_t | \mathfrak{F}_{t-} \right) = \E \left( d \left(\sum_i \frac{x_t^2(i)}{N} \right) | \mathfrak{F}_{t-}. \right)
$$
We have from \protect \MakeUppercase {P}roposition\nobreakspace \ref {M2 Dynamics} that
\begin{align*}
\E ( x_t^2(i) &\vert i \text{ and } j \text{ meet at } t , \mathfrak{F}_{t-} ) \\
&= (1 - \eta)^2 x_{t-}^2(i) + 2 \eta (1 - \eta) x_{t-}(i) x_{t-}(j) + \eta^2 x_{t-}(j).
\end{align*}
When agents $i$ and $j$ meet, only $x(i)$ and $x(j)$ change and so
\begin{align*}
\E ( Q_t - Q_{t-} &\vert i \text{ and } j \text{ meet at } t , \mathfrak{F}_{t-} ) \\
&= - \E \left( \frac{x_t^2(i) - x_{t-}^2(i)}{N} + \frac{x_t^2(j) - x_{t-}^2(j)}{N} \vert i \text{ and } j \text{ meet at } t , \mathfrak{F}_{t-} \right) \\
&= (2\eta - \eta^2) \frac{ x_{t-}^2(i) + x_{t-}^2(j)}{N} - 4 \eta (1 - \eta) \frac{x_{t-}(i) x_{t-}(j)}{N} \\
&- \eta^2 \frac{ x_{t-}(j) + x_{t-}(i)}{N} \\
&= \frac{4 \eta( 1 - \eta)}{2N}(x_{t-}(i) - x_{t-}(j))^2 \\
&- \frac{\eta^2}{N} ( x_{t-}(i) (1 - x_{t-}(i)) + x_{t-}(j)(1 - x_{t-}(j)).
\end{align*}
Summing over $i$ and $j$ our first equation for $d Q_t$ is done. The second is an immediate consequence of the first using the identity
$$
\varepsilon(x, x)_t \geq \lambda \var(x)_t = \lambda(M_t(1 - M_t) - Q_t).
$$
\end{proof}
\end{prop}

\subsection{Escaping an $\epsilon$ Neighbourhood}

Next we focus our attention on how long it takes $M_t$ to escape from the neighbourhood $\left( M_0 - \epsilon, M_0 + \epsilon \right)$ for some small (unspecified for now) $\epsilon$. Let $\tau$ be the escape time, i.e.
$$
\tau = \inf \{ t \geq 0 \colon M_t \notin \left( M_0 - \epsilon, M_0 + \epsilon \right) \}.
$$

For ease of notation in this section we will often write $\E$ for $\E_{x_0}$ - that is the expectation starting from some  initial condition $x_0$, perhaps with some (to be specified) condition on $M_0$.

Our main goal in this section will be to show the following bound.

\begin{prop}{
There exists a positive constant $A(\eta)$ so that for any $M_0$ and $\epsilon$ satisfying
$$
\frac{\eta}{2N} \leq \epsilon \leq \frac{1}{2} M_0 (1 - M_0)
$$
the first exit time $\tau$ satisfies
$$
\E \tau \leq  A(\eta)\frac{N}{\lambda} M_0 (1 - M_0).
$$}
\label{epsilon Neighbourhood Bound}
\end{prop}

\subsubsection{Lower Bound for $\E M_{\tau}^2 - M_0^2$}
First we look for a bound on the heterozygosity $M_t(1 - M_t)$. We will make repeated use of the following calculus exercise.

\begin{lemma}{For a fixed $x_0$, if
$$
\epsilon \leq \frac{x_0 (1 - x_0)}{2}
$$
and $x_0 - \epsilon \leq x \leq x_0 + \epsilon$ then
$$
x(1 - x) \geq \frac{1}{2}x_0 (1 - x_0).
$$}
\label{M(1 - M) Bound}
\end{lemma}

Using our process dynamics calculations we may now begin to bound $\tau$.

\begin{lemma}{There exist positive constants $C(\eta), D(\eta)$ so that
$$
\E  \int_0^{\tau} Q_r dr \geq C(\eta)  \lambda M_0(1 - M_0) \E \tau - D(\eta)(\E Q_{\tau} - Q_0).
$$}
\label{Q Integral Bound}
\begin{proof}
First we recall that from \protect \MakeUppercase {P}roposition\nobreakspace \ref {Q Dynamics} we have a submartingale
$$
Y_t = Q_t - Q_0 - 4 \lambda \eta(1 - \eta) \int_0^t M_r (1 - M_r) dr + (2 \eta^2 + 4 \lambda \eta(1 - \eta)) \int_0^t Q_r dr.
$$
The Optional Stopping Theorem  shows $\E Y_{\tau} \geq \E Y_0 = 0$, so
\begin{align*}
\E Q_{\tau} - &Q_0 + \left( 2 \eta^2 + 4 \lambda \eta(1 - \eta) \right) \int_0^{\tau} Q_r dr  \\
&\geq  4 \lambda \eta(1 - \eta) \E \int_0^{\tau} M_r (1 - M_r) dr   \\
&\geq  4 \lambda \eta(1 - \eta) \E \int_0^{\tau} \frac{1}{2}M_0 (1 - M_0) dr  \text{ by Lemma\nobreakspace \ref {M(1 - M) Bound}}\\
&\geq 2 \lambda \eta( 1- \eta) M_0 (1 - M_0) \E \tau.
\end{align*}

Next, we note that since $\lambda \leq 1$
$$
2 \eta^2 + 4 \lambda \eta (1 - \eta) \leq 4 \eta - 2 \eta^2.
$$
Substituting this in and rearranging the inequality 
$$
(4 \eta - 2 \eta^2) \E \int_0^{\tau} Q_r dr \geq 2 \lambda \eta(1 - \eta) M_0(1 - M_0) \E \tau - \E Q_{\tau} + Q_0
$$
and so
$$
 \E \int_0^{\tau} Q_r dr \geq C(\eta) \lambda M_0(1- M_0)\E \tau - D(\eta)( \E Q_{\tau} - Q_0 )
$$
for $C(\eta) = \frac{2 \eta(1 - \eta)}{4 \eta - 2 \eta^2}$ and $D(\eta) = \frac{1}{4 \eta - 2 \eta^2}$.
\end{proof}
\end{lemma}

Using this we are ready for our lower bound.

\begin{lemma}{There exist positive constants $A(\eta), B(\eta)$ with 
$$
\E M_{\tau}^2 - M_0^2 \geq \frac{1}{N}\left( A(\eta)\lambda  M_0(1 - M_0) \E \tau - B(\eta)\left( \E Q_{\tau} - \E Q_0 \right) \right).
$$}
\label{M2 Lower Bound}
\begin{proof}
\protect \MakeUppercase {P}roposition\nobreakspace \ref {M2 Dynamics} shows
$$M_t^2 - M_0^2 - \frac{2 \eta^2}{N}\int_0^t Q_r dr$$
is a martingale. The Optional Stopping Theorem and Lemma\nobreakspace \ref {Q Integral Bound} show
$$
\E M_{\tau}^2 - M_0^2 \geq \frac{2 \eta^2}{N}
 \left( C(\eta) \lambda M_0(1- M_0)\E \tau - D(\eta)( \E Q_{\tau} - Q_0 ) \right),
$$
which finishes our proof.
\end{proof}
\end{lemma}

\subsubsection{Upper Bound for $\E M_{\tau}^2 - M_0^2$}

We will now impose the condition on our $\epsilon$ neighbourhoods that
$$
\frac{\Delta}{4} \leq \epsilon \leq \frac{M_0(1 - M_0)}{2},
$$
where $\Delta = \frac{2 \eta}{N}$ is the maximum step size for $M_t$. 

\begin{lemma}{For $\tau$ the first escape time from $(M_0 - \epsilon, M_0 + \epsilon)$ we have
$$
\E M_{\tau}^2 - M_0^2 \leq \frac{25}{4} M_0^2 (1 - M_0)^2
$$}
\label{M2 Upper Bound}
\begin{proof}
Consider the martingale $Y_t = M_t - M_0$. From the bound on step sizes, we have that at $\tau$ 
$$
|Y_{\tau}| \leq \epsilon + \Delta \leq 5 \epsilon
$$
and so
$$
\E Y_{\tau}^2 \leq 25 \epsilon^2 \leq \frac{25}{4} M_0^2 (1 - M_0)^2.
$$
From the Optional Stopping Theorem
$$
\E Y_{\tau}^2 = \E (M_{\tau} - M_0)^2 = \E M_{\tau}^2 - M_0^2.
$$

\end{proof}
\end{lemma}

\subsubsection{Bounding $\E \tau$}
We can now prove \protect \MakeUppercase {P}roposition\nobreakspace \ref {epsilon Neighbourhood Bound}.

\begin{proof}
Combining Lemma\nobreakspace \ref {M2 Lower Bound} and Lemma\nobreakspace \ref {M2 Upper Bound} in opposite ways, we have
$$
\frac{25}{4} M_0^2(1 - M_0)^2 \geq \frac{1}{N}\left( A(\eta)\lambda M_0(1 - M_0)  \E \tau - B(\eta) \left( \E Q_{\tau} - \E Q_{0} \right) \right).
$$
Rearranging this and redefining the constants $A(\eta)$ and $B(\eta)$ we have that
$$
\E \tau \leq A(\eta) \frac{N}{\lambda} M_0 (1 - M_0) + \frac{B(\eta)}{\lambda M_0 (1 - M_0)}\left( \E Q_{\tau} - \E Q_{0} \right).
$$
Noting that by assumption $M_0(1 - M_0) \geq \frac{ \eta}{N}$ and so $\frac{1}{M_0(1 - M_0)} \leq \frac{N}{\eta^2}$, we find
$$
\E \tau \leq \frac{N}{\lambda} \left( A(\eta) M_0 (1 - M_0) + \frac{B(\eta)}{\eta^2} \left( \E Q_{\tau} - \E Q_{0} \right) \right).
$$
Finally, since $M_t (1 - M_t)$ is a supermartingale the Optional Stopping Theorem gives $\E M_{\tau}(1 - M_{\tau}) \leq M_0(1 - M_0)$. Jensen's inequality implies $Q_t \leq M_t(1 - M_t)$ always and so combining these two facts
$$
\E Q_{\tau} - Q_0 \leq \E M_{\tau}(1 - M_{\tau}) \leq M_0 (1 - M_0),
$$
giving that
$$
\E \tau \leq \frac{N}{\lambda}M_0 (1 - M_0) \left( A(\eta) + \frac{B(\eta)}{\eta^2} \right)
$$
which completes our proof.
\end{proof}

\subsection{Embedding in the Wright-Fisher Diffusion}

We will now analyse the martingale $M_t$, starting from some initial configuration $x_0$, by discretizing it and embedding that into a Brownian motion. We begin by defining a series of stopping times $\tau_k$ for the martingale $M_t$. Let $\tau_0 = 0$ and for $k \geq 1$ define $\tau_k$ inductively as
$$
\tau_k = \inf \left\lbrace t \geq \tau_{k - 1} \colon \vert M_t - M_{\tau_{k - 1}} \vert \geq \frac{M_{\tau_{k - 1}}(1 - M_{\tau_{k - 1}})}{2} \right\rbrace,
$$
that is the first time after $\tau_{k - 1}$ that $M_t$ exits the ball of size $\frac{M_{\tau_{k - 1}}(1 - M_{\tau_{k - 1}})}{2}$ around $M_{\tau_{k - 1}}$.

Using \protect \MakeUppercase {P}roposition\nobreakspace \ref {epsilon Neighbourhood Bound}, we have the following bound on the expectation of the increments of our stopping times.

\begin{lemma}{There exists a constant $A(\eta)$ so that, from any initial $x_0$, assuming that $M_{\tau_{k - 1}} \in ( \Delta, 1 -  \Delta)$:
$$
\E \left( \tau_k - \tau_{k - 1} \vert F_{\tau_{k - 1}} \right)
 \leq A(\eta) \frac{N}{\lambda}  M_{\tau_{k - 1}}(1 - M_{\tau_{k - 1}}).
$$}
\label{tau Increment Bound}
\begin{proof}
This follows immediately from \protect \MakeUppercase {P}roposition\nobreakspace \ref {epsilon Neighbourhood Bound}, applying the Strong Markov property at time $\tau_{k - 1}$.
\end{proof}
\end{lemma}

\subsubsection{The Wright-Fisher Diffusion}
Here we introduce the Wright-Fisher diffusion $W_t$. We recall that $W_t$ is the continuous martingale in $[0, 1]$ solving the stochastic differential equation
$$
 d W_t = \sqrt{W_t (1 - W_t)} d B_t.
$$
We will frequently think of $W_t$ as a diffusion process and so record here that $W_t$ has instantaneous drift $\mu(x) = 0$ and variance $\sigma^2(x) = x (1 - x)$. 

\begin{lemma}{From any initial $W_0 = w_0$, if $\epsilon \leq \frac{w_0(1 - w_0)}{2}$ then the first escape $\tau$ from the $\epsilon$-ball about $w_0$ satisfies
$$
\frac{1}{3} w_0(1 - w_0) \leq \E_{w_0} \tau \leq \frac{5}{3} w_0 (1 - w_0).
$$
}
\label{WF Escape Bound}
\begin{proof}
Let $w_{\pm} = w_0 \pm \epsilon$. On the interval $[w_-, w_+]$ let
$$
u(x) = \E_x \tau.
$$
Applying a standard argument for diffusion processes we find that $u(x)$ satisfies the equation
$$
-1 = \frac{x(1 - x)}{2} u''(x),
$$
subject to the boundary conditions $u(w_-) = u(w_+) = 0$. Integrating, we find that for $$f(x) = 2 \left( x \ln x + (1 - x) \ln(1 - x) \right)$$ we have the solution
$$
u(x) = - f(x) + A x + B
$$
for some constants $A, B$. Applying the boundary conditions we get
$$
2 \epsilon A = f(w_+) - f(w_-) 
$$
and
$$
B = f(w_+) - A w_+.
$$
Therefore, for $w_0$ we have
\begin{align*}
u(w_0) &= - f(w_0) + A w_0 + f(w_+) - A w_+ \\
&= - f(w_0) - A \epsilon + f(w_+) \\
&= - f(w_0) - \frac{f(w_+) - f(w_-)}{2} + f(w_+) \\
&= - f(w_0) + \frac{1}{2}( f(w_+) + f(w_-)). 
\end{align*}

To simplify this we apply the Taylor approximation to $f$ at $w_0$
$$
f(w_0  + \Delta) \approx f(w_0) + f'(w_0) \Delta + \frac{1}{w_0(1 - w_0)} \Delta^2,
$$
which plugging into $u(w_0)$ we find
$$
u(w_0) \approx \frac{\epsilon^2}{w_0(1 - w_0)} 
$$
To approximate the error, we use Taylor's remainder theorem. First,
$$
f^{(3)}(x) = \frac{2}{(1 - x)^2} - \frac{2}{x^2} = \frac{4x - 2}{x^2(1 -x)^2}
$$
and so
$$
| f^{(3)}(x) | \leq \frac{2}{x^2 (1 - x)^2}.
$$
On the interval $[w_0 - \epsilon, w_0 + \epsilon]$, applying a bound we used for $M_t(1 - M_t)$ above, we have
$$
x(1 - x) \geq \frac{w_0(1 - w_0)}{2}.
$$
Therefore the third derivative is bounded on the same interval by
$$
|f^{(3)}(x) | \leq \frac{8}{w_0^2(1 - w_0)^2}.
$$
Therefore the error $R(w_0 + x)$ for $|x| \leq \epsilon$ of the 2nd Taylor approximation is bounded by
$$
|R(w_0 + x)| \leq \frac{8}{w_0^2 (1 - w_0)^2} \frac{|x|^3}{6}.
$$
So, for $w_{\pm}$, using that $\epsilon \leq \frac{w_0(1 - w_0)}{2}$ we have
$$
|R(w_{\pm})| \leq \frac{2 \epsilon^2}{3 w_0 (1 - w_0)}.
$$
Applying this to our calculation of $u(w_0)$ - once to each of $f(w_0 \pm \epsilon)$ we have
\begin{align*}
|u(w_0) - \frac{\epsilon^2}{w_0(1 - w_0)}| &\leq \frac{1}{2}|R(w_0 - \epsilon) + R(w_0 + \epsilon) | \\
&\leq \frac{2 \epsilon^2}{3 w_0 (1 - w_0)}, \\
\end{align*}
from which our bounds on $u(w_0)$ follow easily.
\end{proof}
\end{lemma}

We will also need the well known bound for the absorption time $T_{\abs}$ of $W_t$ at the boundary $\{ 0, 1 \}$. To state this bound, recall the function
\begin{align}
\phi(x) = -x \ln(x) - (1 - x) \ln(1 - x) \label{Define phi}.
\end{align}

\begin{lemma}{Let $T_{\abs}$ be the stopping time
$$
T_{\abs} = \inf \{ t \geq 0 \colon W_t \notin (0, 1) \},
$$
i.e. the time when $W_t$ is absorbing in $\{0, 1\}$. Then $\E_x T_{\abs} = 2 \phi(x) $ for any $x \in [0, 1]$.}
\label{WF Absorption Time}
\begin{proof}
This is easy using the same techniques as the previous proof. Letting
$$
u(x) = \E_x T_{\abs}
$$
we find that
$$
u(x) = - 2 ( x \ln x + (1 - x) \ln(1 - x)) + Ax + B,
$$
subject to the boundary conditions $u(0)= u(1) = 0$. This necessitates $A = B = 0$ and so
$$
u(x) = - 2 ( x \ln x + (1 - x) \ln(1 - x)),
$$
completing the proof.
\end{proof}
\end{lemma}

\subsubsection{Embedding $M_{\tau_k}$ in $W_t$}

 Let $W_t$ be a Wright-Fisher diffusion started at $W_0 = M_0$. As $|M_t| \leq 1$ the discrete martingale $\{ M_{\tau_k} \}_{k \geq 0}$ is clearly square integrable and so \cite{dubins1968} we can find a sequence of stopping times $\tilde{\tau}_k$ for $W_t$ so that
 \begin{equation}
 \{ M_{\tau_k} \}_{k \geq 0} =^d \{ W_{\tilde{\tau}_k} \}_{k \geq 0}. \label{Equivalence in Distribution}
  \end{equation}
 
 We will use this embedding to bound the first escape time $S$ by comparison with the absorption time of the Wright-Fisher diffusion.
 
\subsubsection{The Comparison Calculation} 
 
 We will focus on the first time that $M_t$ exits the interval $(\frac{\eta}{2N}, 1 - \frac{\eta}{2N})$. Recall the stopping time $S$ defined by
 \begin{equation}
 S = \inf \{ t \geq 0 \colon M_t \notin (\frac{\eta}{2N}, 1 - \frac{\eta}{2N}) \}, \label{Define S}
  \end{equation}
 and let 
 $$
 K = \inf \{ k \geq 0 \colon \tau_k \geq S \}.
 $$
 
 \begin{lemma}{ $K < \infty$ almost surely.}
 \label{K Finite}
 \begin{proof}
 For any $0 \leq x \leq \frac{1}{2}$ we have
 $$
 x - \frac{x(1 - x)}{2} \leq \frac{3}{4}x,
 $$
 and so if $M_0 \leq \frac{1}{2}$, then for any $k \geq \frac{\ln(\frac{\eta}{N})}{\ln(3/4)}$ we have $M_{\tau_k} \leq \frac{\eta}{2N}$ with positive probability. Fix such a $k_0$. In fact, as $\{ M_{\tau_n} \}_{n \geq 0}$ is a martingale, we have
 $$
 \Prb \left( M_{\tau_{k_0}} \leq \frac{\eta}{2N} \right) \geq \frac{1}{2^{k_0}}.
 $$
 A similar argument holds for an initial configuration with $\frac{1}{2} \leq M_0 \leq 1$ going above the level $1 - \frac{\eta}{2N}$. Therefore, for any initial $M_0$, there is a uniform lower bound on the probability that $K < k_0$. Thus, by a standard Strong Markov argument, we must have $K < \infty$ almost surely.
 \end{proof}
 \end{lemma}
  
 Next, we define equivalent stopping times for $W_t$. Let $\tilde{K}$ be the index
 $$
 \tilde{K} = \inf \{ k \geq 0 \colon W_{\tilde{\tau}_k} \notin (\frac{\eta}{2N}, 1 - \frac{\eta}{2N})  \}
 $$
 and recall the absorption time $T_{\abs}$ given by
 $$
 T_{\abs} = \inf \{ t \geq 0 \colon W_t \notin (0, 1) \}.
 $$
 For $M_t$ we have clearly that 
 \begin{equation}
 S \leq \tau_K \label{S leq tau_K}.
 \end{equation}
 
  Furthermore, as $M_{\tau_K} \in (0, 1)$ so must be $W_{\tilde{\tau}_{\tilde{K}}}$ by the equivalence in distribution (and thus support). Therefore,
 $$
 \tilde{\tau}_{\tilde{K}} \leq T_{\abs},
 $$
 as for  $t \geq T_{\abs}$, $W_t$ is constant and in $\{ 0, 1 \}$.

	First we need to bound the $\tilde{\tau}_k$.
	
	\begin{lemma}{The hitting times $\tilde{\tau}_k$ satisfies
	$$
	\E \left( \tilde{\tau}_k - \tilde{\tau}_{k - 1} \vert F_{\tilde{\tau}_{k -1}} \right) \geq \frac{1}{3} W_{\tilde{\tau}_{k -1}} (1 - W_{\tilde{\tau}_{k -1}}).
	$$}
	\begin{proof}
	For this, we note that starting at $w_0 = W_{\tilde{\tau}_{k - 1}}$, the time $\tilde{\tau}_k$ can only occur after $W_t$ leaves the interval $(W_{\tilde{\tau}_k} - \frac{w_0(1 - w_0)}{2}, W_{\tilde{\tau}_k} + \frac{w_0(1 - w_0)}{2} )$ as $W_{\tilde{\tau}_k}$ is already outside this interval and $W_t$ is continuous. Write $\tau$ for the first exit time of this interval. Applying the Strong Markov Property we see that
	\begin{align*}
	\E \left( \tilde{\tau}_k - \tilde{\tau}_{k - 1} \vert F_{\tilde{\tau}_{k -1}} \right) 
	&\geq \E \left( \tau \vert F_{\tilde{\tau}_{k - 1}} \right) \\
	&\geq \frac{1}{3} W_{\tilde{\tau}_{k -1}} (1 - W_{\tilde{\tau}_{k -1}}) \text{ by Lemma\nobreakspace \ref {WF Escape Bound}},
	\end{align*}
	completing our proof.
\end{proof} 
\end{lemma}

We are now ready to prove $\protect \MakeUppercase {T}heorem\nobreakspace \ref {Main Escape Theorem}$.
 
 \begin{proof}
 We recall by Equation\nobreakspace \textup {(\ref {S leq tau_K})}, $\E S \leq \E \tau_K$ and so we will focus on bounding $\E \tau_K$. As $K < \infty$ almost surely by Lemma\nobreakspace \ref {K Finite} we have that
 $$
\E \tau_K = \E \sum_{k = 1}^{\infty} (\tau_{k} - \tau_{k - 1})1_{K \geq k}. 
 $$
 By Lemma\nobreakspace \ref {tau Increment Bound}
 $$
 \E \left( \tau_k - \tau_{k - 1} \vert F_{\tau_{k - 1}} \right) 
 \leq A(\eta) \frac{N}{\lambda} M_{\tau_{k - 1}} (1 - M_{\tau_{k - 1}}),
 $$
 for some constant $A(\eta)$ depending only on $\eta$. Therefore we can calculate using the Strong Markov property that
 \begin{align*}
 \E (\tau_{k} - \tau_{k - 1})1_{K \geq k} &= \E \E \left( (\tau_{k} - \tau_{k - 1})1_{K \geq k} \vert F_{\tau_{k - 1}}\right) \\
 &= \E 1_{K \geq k} \E_{x_{\tau_{k - 1}}} \left( \tau_{k} - \tau_{k - 1} \right) \\
 &\leq  A(\eta) \frac{N}{\lambda} \E 1_{K \geq k} M_{\tau_{k -1}} (1 - M_{\tau_{k - 1}}) .
 \end{align*}

	From Equation\nobreakspace \textup {(\ref {Equivalence in Distribution})} $\{ M_{\tau_k} \}_{k \geq 0}$ and $\{ W_{\tilde{\tau}_k} \}_{k \geq 0}$ are equivalent in distribution, so
	\begin{align*}
	\E 1_{K \geq k} M_{\tau_{k -1}} (1 - M_{\tau_{k - 1}}) &=
	 \E 1_{\tilde{K} \geq k} W_{\tilde{\tau}_{k - 1}} (1 - W_{\tilde{\tau}_{k - 1}} ).
	 \end{align*}
	 
	 By Lemma\nobreakspace \ref {WF Escape Bound} 
	 $$
	\frac{1}{3} W_{\tilde{\tau}_{k - 1}} (1 - W_{\tilde{\tau}_{k - 1}} ) \leq \E \left( \tilde{\tau}_k - \tilde{\tau}_{k - 1} \vert F_{\tilde{\tau}_{k - 1}} \right),
	 $$
	so we can calculate
	 \begin{align*}
	 \E 1_{\tilde{K} \geq k} W_{\tilde{\tau}_{k - 1}} (1 - W_{\tilde{\tau}_{k - 1}} )
	 &\leq \E 1_{\tilde{K} \geq k} 3\E \left(\tilde{\tau}_k - \tilde{\tau}_{k - 1} \vert F_{\tilde{\tau}_{k - 1}}\right) \\
	 &= 3 \E \E \left( (\tilde{\tau}_k - \tilde{\tau}_{k - 1})1_{\tilde{K} \geq k} \vert F_{\tilde{\tau}_{k - 1}}\right) \\
	 &=3 \E (\tilde{\tau}_k - \tilde{\tau}_{k - 1})1_{\tilde{K} \geq k}.
	\end{align*} 

Therefore we see that
\begin{align*}
 \E \tau_S &\leq  3 \frac{N}{\lambda} A(\eta) \E \sum_{k \geq 0} \E (\tilde{\tau}_k - \tilde{\tau}_{k - 1})1_{\tilde{K} \geq k} \\
&\leq 3 \frac{N}{\lambda} A(\eta) \E \tilde{\tau}_{\tilde{K}} \\
&\leq 3 \frac{N}{\lambda} A(\eta) \E T_{\abs}\\
\end{align*}

 Because Equation\nobreakspace \textup {(\ref {S leq tau_K})} $S \leq \tau_K$ and using Lemma\nobreakspace \ref {WF Absorption Time} to bound $\E T_{\abs}$ we can conclude that
\begin{align*}
\E S &\leq   3 \frac{N}{\lambda} A(\eta) \E_{W_0} T_{\abs} \\
&= 6 A(\eta) \frac{N}{\lambda} \phi(W_0) \\
&= 6 A(\eta) \frac{N}{\lambda} \phi(M_0)
\end{align*}
from which our conclusion follows.
 \end{proof}

\section{The General Model}
\label{General Model Section}

We are now ready to prove our bound on the fixation time of the general iPod model with $\sigma$ songs. We recall that for each agent $i$, we write their preference for song $k$ by $X_t^k(i)$. For each song $k$, we write the average preference for that song as $M_t^k$, given by
$$
M_t^k = \sum_i \frac{X_t^k(i)}{N}.
$$
We have shown in \protect \MakeUppercase {P}roposition\nobreakspace \ref {M Martingale} that for each $1 \leq k \leq \sigma$, $M_t^k$ is a martingale.

\subsection{Approaching the Boundary}

We begin by using \protect \MakeUppercase {T}heorem\nobreakspace \ref {Main Escape Theorem} to bound the time it takes for one of the $M^k_t$ to approach the boundary $1$. Specifically, we will analyse the stopping time
\begin{align}
\label{Define tau}
\tau = \inf \{t \geq \colon M^k_t \geq 1 - \frac{\eta}{2 N} \text{ for some } 1 \leq k \leq \sigma \}.
\end{align}

Let $S^k$ be the first time that $M^k_t$ approaches either boundary, that is
$$
S^k = \inf \{ t \geq 0 \colon M^k_t \notin (\frac{\eta}{2N}, 1 - \frac{\eta}{2N}) \},
$$
and set $S_{\max} = \max_{1 \leq k \leq \sigma} S^k$. We will bound $\tau$ in two steps, first by bounding $\E S_{\max}$ and second by showing that $\E \tau$ is on the same order of magnitude as $\E S_{\max}$.

\begin{prop}{There is a constant $C(\eta)$ so that from any initial configuration $X_0$ we have
$$
\E S_{\max} \leq C(\eta) \ln(\sigma) \frac{N}{\lambda}.
$$}
\label{Smax bound}
\begin{proof}
To begin we recall that by \protect \MakeUppercase {T}heorem\nobreakspace \ref {Main Escape Theorem} there is a constant $C(\eta)$ so that
$$
\E S^k \leq C(\eta) \frac{N}{\lambda} \phi(M_0^k),
$$
for all $1 \leq k \leq \sigma$. Clearly $S_{\max} \leq \sum_{k = 1}^{\sigma} S^k$ and so
$$
\E S_{\max} \leq C(\eta) \frac{N}{\lambda} \sum_{k =1 }^{\sigma} \phi(M_0^k).
$$
Recalling that the $M_0^k$ satisfy the constraint $\sum_{k = 1}^{\sigma} \phi(M_0^k) = 1$, a simple calculus exercise in Lagrange multipliers shows that $\sum_{k =1 }^{\sigma} \phi(M_0^k)$ is maximized when $M_0^k = \frac{1}{\sigma}$ for all $k$ and so
\begin{align*}
\sum_{k =1 }^{\sigma} \phi(M_0^k) &= \sigma \left(- \frac{1}{\sigma} \ln(\frac{1}{\sigma}) - (1 - \frac{1}{\sigma}) \ln(1 - \frac{1}{\sigma}) \right) \\
&\leq \ln(\sigma) + 1,
\end{align*}
completing the proof.
\end{proof}
\end{prop}

Next we will show that $\E \tau$ is on the same order of magnitude (w.r.t $N$) as $\E S_{\max}$.

\begin{prop}{From any initial configuration $X_0$
$$
\E \tau \leq 2 \sup_{X_0} \E_{X_0} S_{\max}.
$$}
\label{tau Bound}
\begin{proof}
For $1 \leq k \leq \sigma$ let 
$$
A_k = \{ M^k_{S_k} \geq 1 - \frac{\eta}{2N} \},
$$
that is the event that song $k$ approaches the boundary $1$ at time $S_k$ (as opposed to the boundary $0$). Clearly $A_k \subset \{ \tau \leq S_{\max} \}$ and so
\begin{align*}
\Prb( \tau \leq S_{\max} ) & \geq \Prb( \cup_{1 \leq k \leq \sigma} A_k ).
\end{align*}

We claim that the $A_k$ are almost surely disjoint. Assuming otherwise, given $A_k \cap A_j$ one of $S^k, S^j$ must occur first - since both $M^k$ and $M^j$ can't be greater than $1 - \frac{\eta}{2N}$ at the same time. Assume without loss of generality that $S^k < S^j$. At time $S^k$, $M^k \geq 1 - \frac{\eta}{2N}$ and so we must have $M^j \leq \frac{\eta}{2N}$, meaning that $S^j$ has already occurred, a contradiction. Therefore the events $A_k$, $1 \leq k \leq \sigma$ are almost surely disjoint.

Next, as $M^k_t$ is a martingale - with step size bounded by $\Delta = \frac{2 \eta}{N}$ - we have
$$
\Prb(A_k) \geq \frac{M^k_0 - \frac{\eta}{2N}}{1 - \frac{\eta}{N} + \Delta},
$$
and so for $N >> 0$
$$
\Prb(A_k) \geq \frac{M^k_0}{2}.
$$

 Thus for any initial configuration $X_0$
$$
\Prb_{X_0}(\tau \leq S_{\max}) = \sum_k \Prb(A_k) \geq \sum_k \frac{	M^k_0}{2} = \frac{1}{2},
$$
or equivalently $\Prb_{X_0}(\tau \geq S_{\max}) \leq \frac{1}{2}$.

Let $m = \sup_{X_0} \E_{X_0} \tau$. Applying the Strong Markov Property we have
\begin{align*}
\E_{X_0} \tau &\leq \E_{X_0} S_{\max} + \E_{X_0} (\tau - S_{\max})1_{\tau \geq S_{\max}} \\
&\leq \E_{X_0} S_{\max} + \E 1_{\tau \geq S_{\max}} \E_{X_{S_{\max}}} \tau \\
&\leq \E_{X_0} S_{\max} + \Prb_{X_{S_max}}( \tau \geq S_{\max} ) m \\
&\leq \E_{X_0} S_{\max} + \frac{m}{2}
\end{align*}
which implies that
$$
m \leq 2 \sup_{X_0} \E_{X_0} S_{\max}.
$$
\end{proof}
\end{prop}

\subsection{Fixation Time}

Next we will estimate the fixation time given that the preference $M_t^k$ for some (fixed but arbitrary) song $k$ has approached the boundary $1$. Specifically, we will consider starting from an initial configuration $X_0$ with
$$
M^k_0 \geq 1 - \frac{\eta}{2N}.
$$
When $M^k$ is near $1$, the fixation time $T_{\fix}$ can only be the last time any song other than $k$ plays. Of course this need not occur. Projecting on $k$, this is the last time one of the Bernoulli trials for $k$ has failed. We begin by showing that from such an initial configuration, $T_{\fix}$ has likely already occurred. 

\begin{prop}{From an initial configuration $X_0$ with $M^k_0 \geq 1 - \frac{\eta}{2N}$, we have
$$
\Prb_{X_0} \left( T_{\fix} = 0 \right) \geq \frac{1}{2}
$$}
\label{Tfix Probability}
\begin{proof}
 We will consider the stopping time $R$, the first time any song other than $k$ plays. Before $R$, each $X^k(i)$ can only increase. Therefore at time $R$ - without loss of generality, a meeting of $i$ and $j$ - if another song is played by only one of $i$, $j$ then
\begin{align}
X_R^k(i) + x_R^k(j) &= (1 - \eta)(X_{R-}(i) + X_{R-}(j)) + \eta \\
&\leq 2(1 - \eta) + \eta \\
&= 2 - \eta.
\end{align}
If both agents play a different song, then $X^k(i) + X^k(j)$ is even smaller at $R$.

This then implies that on $\{ R < \infty \}$
$$
M^k_R \leq 1 - \frac{\eta}{N}.
$$

Now, applying the Optional Stopping Theorem to $R \wedge t$, we find that
\begin{align}
1 - \frac{\eta}{2N} & \leq M_0 \\
&= \E M^k_{R \wedge t} \\
&= \E \left( M^k_R 1_{R \leq t} + M^k_t 1_{t < R} \right) \\
&\leq (1 - \frac{\eta}{N})(1 - \Prb \left( t < R \right)) + 1 \Prb \left( t < R \right).
\end{align}
Solving for $\Prb \left(t < R \right)$ we find that
$$
\Prb\left( t < R \right) \geq \frac{1}{2}.
$$
As this is true for arbitrary $t$, we have $\Prb \left( R = \infty \right) \geq \frac{1}{2}$ from which our result follows.
\end{proof}
\end{prop}

Next we need to consider what happens when $M_t^k$ approaches $1$, but the song $k$ fails to play at a meeting. 

\begin{prop}{Consider the stopping time $R$ given by
$$
R = \inf \{t \geq 0 \colon \text{ some song other than k plays at } t  \}.
$$
From any initial configuration $M^k_0 \geq 1 - \frac{\eta}{2 N}$, we have
$$
\E R 1_{R < \infty} \leq \frac{1}{8 \eta}.
$$}
\begin{proof}
Let $T_n$, $1 \leq n < \infty$ be the $n$-th meeting time. We first define
$$
\tilde{R} = \inf \{n \geq 0 \colon \text{ some song other than k plays at } T_n  \}.
$$
We will calculate how $M^k_t$ changes after the first meeting time, given that song $k$ is played by both agents at the meeting time $T_1$.

If agents $i$ and $j$ meet and both play $k$ at $T_1$ then
$$
X^k_{T_1}(i) = (1 - \eta) X^k_0(i) + \eta
$$
and similarly for $X^k(j)$. So given that $i$ and $j$ meet and play $k$
$$
M^k_{T_1} = M^k_0 - \frac{\eta(X^k_0(i) + X^k_0(j))}{N} + \frac{2 \eta}{N}.
$$
Summing over pairs of agents we find that
\begin{align*}
\E & \left(M^k_{T_1}  \vert \text{ both agents play k at } T_1 , \mathfrak{F}_0 \right) \\ 
&= \sum_{i, j} \E \left( M^k_{T_1} \vert \text{ i meets j, both play k at } T_1, \mathfrak{F}_0 \right)
\Prb \left(\text{ i meets j at } T_1 \vert \mathfrak{F}_0 \right) \\
&= \sum_{i, j} \frac{\nu_{ij}}{N}  \E \left( M^k_0 - \frac{\eta(X^k_0(i) + X^k_0(j) - 2)}{N} \vert \text{ i \& j both play k at } T_1, \mathfrak{F}_0 \right) \\
&= \sum_{i, j} \frac{\nu_{ij}}{N} (M^k_0 - \frac{\eta(X^k_0(i) + X^k_0(j) - 2)}{N}) \\
&= M^k_0 + \frac{2 \eta}{N} - \sum_{i,j} \frac{\nu_{ij}}{N} \frac{\eta(X^k_0(i) + X^k_0(j))}{N} \\
&= M^k_0 + \frac{2 \eta}{N} - \frac{2 \eta M^k_0}{N} \\
&= (1 - \frac{2 \eta}{N})M^k_0 + \frac{2 \eta}{N}.
\end{align*}

By the same calculation we find that
$$
\E \left(M^k_{T_2}  \vert \text{both agents play k at } T_2 , \mathfrak{F}_{T_1} \right) = (1 - \frac{2 \eta}{N})M^k_{T_1} + \frac{2 \eta}{N}
$$
and so
\begin{align*}
\E  & \left(M^k_{T_2}  \vert \text{ both agents play k at } T_1 \text{ and } T_2 , \mathfrak{F}_0 \right) \\
&= (1 - \frac{2 \eta}{N})\left((1 - \frac{2 \eta}{N})M^k_0 + \frac{2 \eta}{N} \right) + \frac{2 \eta}{N} \\
&= (1 - \frac{2 \eta}{N})^2 M^k_0 + 1 - (1 - \frac{2 \eta}{N})^2 \\
&= 1 - (1 - \frac{2 \eta}{N})^2 (1 - M^k_0).
\end{align*}

Continuing the same easy inductive calculation we find that
$$
\E  \left(M^k_{T_n}  \vert \tilde{S} > n , \mathfrak{F}_0 \right) = 1 - (1 - \frac{2 \eta}{N})^n (1 - M^k_0).
$$

Next, we need to know the chance of some song other than $k$ being played at time $T_n$ given $M^k_{T_{n - 1}}$. We will need the identity
$$
1 - xy \leq (1 - x) + (1 - y)
$$
for $x, y \leq 1$ - which follows easily from $1 + (1 - x)(1 - y) \geq 1$. Using that, and that the probability of at least one of $i, j$ not playing $k$ is $1 - X^k(i) X^k(j)$,  we have
\begin{align*}
\Prb & \left(\text{ A song other than k is played at } T_n \vert M^k_{T_{n - 1}} \right) \\
&= \sum_{i, j} \frac{\nu_{ij}}{N} \Prb\left(\text{ Another song is played at } T_n \vert M^k_{T_{n - 1}}, \text{i meets j at } T_n \right) \\
&= \sum_{i, j} \frac{\nu_{ij}}{N} \left( 1 - X^k_{T_{n - 1}}(i)X^k_{T_{n - 1}}(j) \right)\\
&\leq \sum_{i,j} \frac{\nu_{ij}}{N} \left( 1 - X^k_{T_{n - 1}}(i) + 1 - X^k_{T_{n - 1}}(j) \right) \\
&\leq 2(1 - M^k_{T_{n - 1}}).
\end{align*}

Therefore we have that
\begin{align*}
\Prb & \left( \tilde{R} = n \vert \mathfrak{F}_0 \right) \\
&= \Prb \left( \tilde{R} > n - 1, \text{Another song is played at } T_n \vert \mathfrak{F}_0 \right) \\
&\leq \Prb \left( \text{Another song is played at } T_n \vert \tilde{R} > n - 1, \mathfrak{F}_0 \right) \\
&\leq \E \left( 2(1 - M^k_{T_{n - 1}}) \vert \tilde{R} > n - 1, \mathfrak{F}_0 \right) \\
&= 2(1 - \frac{2 \eta}{N})^{n - 1} (1 - M^k_0).
\end{align*}
For the first inequality here we used the simple bound
$$
\Prb \left(A \cap B \right) \leq \Prb \left( A \vert B \right).
$$

This allows us to calculate that
\begin{align*}
\E \left( \tilde{R}1_{\tilde{R} < \infty} \vert \mathfrak{F}_0 \right) &= \sum_{n \geq 0} n \Prb \left(\tilde{R} = n \vert \mathfrak{F}_0 \right) \\
&\leq \sum_{n \geq 0} n 2(1 - \frac{2 \eta}{N})^{n - 1} (1 - M^k_0) \\
&= 2 (1 - M^k_0) \sum_{n \geq 0} n (1 - \frac{2 \eta}{N})^{n - 1} \\
&\leq \frac{ \eta}{2 N} \frac{N^2}{4 \eta^2} \\
&=\frac{N}{8 \eta},
\end{align*}
using our assumption that $M^k_0 \geq 1 - \frac{\eta}{2 N}$ and the Taylor series expansion
$$
\sum_{n \geq 0} n x^{n - 1} = \frac{1}{(1 - x)^2},
$$
for $|x| < 1$.

Our result then follows since meetings occur independently at rate $\frac{1}{N}$ and so
$$
\E \left( R 1_{R < \infty} \vert \mathfrak{F}_0 \right) = \frac{1}{N} \E \left( \tilde{R} 1_{\tilde{R} < \infty} \vert \mathfrak{F}_0 \right) .
$$

\end{proof}
\end{prop}

We are finally prepared to prove \protect \MakeUppercase {T}heorem\nobreakspace \ref {Main Theorem}.

\begin{proof}
We will calculate here an upper bound for
$$
m = \max_{X_0} \E_{X_0} T
$$
i.e. the upper bound over all initial configurations $X_0$.

Let $\tau$ be the stopping time from Equation\nobreakspace \textup {(\ref {Define tau})}, i.e. the first time that some song $k$ has $M_t^k \geq 1 - \frac{\eta}{2 N}$ and let $K$ be that song. Note that this defines $K$ uniquely as $1 - \frac{\eta}{2N} \geq \frac{1}{2}$. Let $R$ be stopping time (as above) defined by
$$
R = \inf \{t \geq \tau \vert \text{ some song other than K is played}\}.
$$

We first recall from \protect \MakeUppercase {P}roposition\nobreakspace \ref {Tfix Probability} that at time $\tau$, we have
$$
\Prb_{X_{\tau}} \left( T_{\fix} = 0 \right) \geq \frac{1}{2}.
$$
Also, at time $\tau$, if $T_{\fix}$ has not yet occurred, then some song other than $K$ will play again and so $R < \infty$.

Combining \protect \MakeUppercase {P}roposition\nobreakspace \ref {Smax bound} and \protect \MakeUppercase {P}roposition\nobreakspace \ref {tau Bound} we have that there exists a constant $C(\eta)$ so that from any initial configuration $X_0$
$$
\E_{X_0} \tau \leq C(\eta) \frac{\ln (\sigma) N}{\lambda}.
$$

We then have for any initial $X_0$:
\begin{align*}
\E_{X_0} T &= \E_{X_0} \E\left( (T_{\fix} - \tau) + \tau \vert \mathfrak{F}_{\tau} \right) \\
&= \E_{X_0} \tau + \E_{X_0} \E_{X_{\tau}} T_{\fix} \\
&= \E_{X_0} \tau + \E_{X_0} \E_{X_{\tau}} T_{\fix} 1_{T_{\fix} > 0} \\
&= \E_{X_0} \tau + \E_{X_0} \E_{X_{\tau}} \left( (T_{\fix} - R) 1_{R < \infty} + R 1_{R < \infty} \right) \\
&= \E_{X_0} \tau + \E \E_{X_{\tau}} R 1_{R < \infty} + \E \E \left( (T_{\fix} - R) 1_{R < \infty} \vert R \right) \\
&= \E_{X_0} \tau + \frac{1}{8 \eta} + \E \left( 1_{R < \infty} \E_{X_R} T \right) \\
&\leq C(\eta)\frac{ \ln(\sigma) N }{\lambda} + \frac{1}{8 \eta} +  \E \left( 1_{R < \infty} m \right) \\
&\leq 2C(\eta)\frac{ \ln(\sigma) N }{\lambda}  + \frac{1}{2} \max_{x_0} \E_{x_0} T_{\fix} .
\end{align*}

Here the $\frac{1}{8 \eta}$ is clearly dominated by the first term. Therefore, we have that
$$
\max_{X_0} \E_{X_0} T_{\fix} \leq 2C(\eta)\frac{ \ln(\sigma) N }{\lambda}  + \frac{1}{2} \max_{X_0} \E_{X_0} T_{\fix}
$$
and so
$$
\E_{X_0} T_{\fix} \leq 4 C(\eta)\frac{ \ln(\sigma) N }{\lambda}
$$
from which our conclusion follows.
\end{proof}

\section{The Interaction Parameter $\eta$}

Our goal here is find the asymptotic of our bound with respect to $\eta$. Tracing through the steps of our proof of \protect \MakeUppercase {T}heorem\nobreakspace \ref {Main Escape Theorem}, we may actually prove the following improved bound.

\begin{prop}{There exists a constant $C$ so that from any initial configuration $x_0$, the first escape time $S$  satisfies
$$
\E_{x_0} S \leq \frac{C}{\eta^3(1 - \eta)} \frac{N}{\lambda}.
$$}
\label{Improved S Bound}
\end{prop}

Then, repeating the arguments in Section\nobreakspace \ref {General Model Section}, we may improve our bound in \protect \MakeUppercase {T}heorem\nobreakspace \ref {Main Theorem} on the expectation of the fixation time $T_{\fix}$.

\begin{thm}{There exists a constant $C$ so that from any initial $X_0$ the fixation time $T_{\fix}$ satisfies
$$
\E T_{\fix} \leq \frac{C}{\eta^3 (1 - \eta)} \frac{\ln(\sigma) N}{\lambda}.
$$}
\label{Improved Main Theorem}
\end{thm}

We conjecture that this can actually be improved to depend on $\eta$ as $\frac{1}{\eta ( 1 - \eta)}$.

\section{The Complete Graph Case}

As an example of a geometry in which more can be said than \protect \MakeUppercase {T}heorem\nobreakspace \ref {Main Theorem}, we look at the complete graph $K_N$ on $N$ vertices. Specifically, we have uniform meeting rates between agents, that is $\nu_{ij} = \frac{1}{N - 1}$ for all pairs of agents $i, j$. It is standard fact that the spectral gap $\lambda_{K_N} = 1$ and so \protect \MakeUppercase {T}heorem\nobreakspace \ref {Main Theorem} shows that the fixation time has 
$$\E T_{\fix} = O(N).$$

A simple argument will show that this order of magnitude bound is in fact tight.

\subsection{A Lower Bound}

 Throughout this section we assume that there are at least two songs, i.e. $\sigma \geq 2$. To achieve any reasonable lower bound, we need to ignore starting conditions that are likely already at fixation by time $t = 0$. We call an initial configuration \textbf{non-trivial} if there exists at least one song $k$ with
 $$
 \frac{1}{2 \sigma} \leq M^k_0 \leq 1 - \frac{1}{2 \sigma},
 $$
 and will consider only non-trivial initial configurations. The choice of the factor of $\frac{1}{2}$ here is of course arbitrary. 

\begin{thm}{There exists a constant $C(\eta, \sigma)$ such that for $K_N$ started from any non-trivial initial configuration, the fixation time $T_{\fix}$ has
$$
\E T_{\fix} \geq C(\eta, \sigma) N.
$$ }
\label{Complete Graph Lower Bound}
\begin{proof}
Recalling \protect \MakeUppercase {T}heorem\nobreakspace \ref {Main Escape Theorem}, first consider any one song and consider its average preference $M_t, t \geq 0$. From the proof of \protect \MakeUppercase {P}roposition\nobreakspace \ref {M2 Dynamics}
$$
\E \left( d M_t (1 - M_t) \vert \mathfrak{F}_{t-} \right) = - \frac{2 \eta^2}{N} Q_t dt,
$$
which combined with $Q \leq \frac{1}{4}$ gives that
$$
M_t (1 - M_t) - M_0(1 - M_0) + \frac{\eta^2}{2 N} t
$$
is a sub-martingale.

By assumption, there exists at least one song $k$ with $M^k_0 (1 - M^k_0) \geq \frac{1}{4 \sigma}$. Let 
$$
T_{2} = \inf_{t \geq 0} \{ M^k_t \notin \left( \frac{1}{8 \sigma}, 1 - \frac{1}{8 \sigma} \right) \},
$$ 
be the first time that $M^k_t$ leaves the interval $\left( \frac{1}{8 \sigma}, 1 - \frac{1}{8 \sigma} \right)$. Then applying the Optional Stopping Theorem
$$
\E M^k_{T_2} (1 - M^k_{T_2}) + \frac{\eta^2}{2 N} \E T_2 \geq M_0 (1 -M_0) \geq \frac{1}{4 \sigma}.
$$
At time $T_2$, we have have
$$
 M^k_{T_2} (1 - M^k_{T_2}) \leq \frac{1}{8 \sigma},
$$
and so we can conclude that
$$
\E T_2 \geq \frac{N}{4  \eta^2 \sigma}.
$$
To complete the proof, we need only show that the fixation time $T_{\fix}$ is with high probability the same order of magnitude as $T_2$.

Consider the first meeting after time $T_2$, between some agents $i$ and $j$. If two different songs are played at that meeting, then by definition $T_{\fix}$ must not have yet occurred. The probability that at a meeting at time $t$ that agent $i$ plays song $k$ and $j$ does not, or vis-versa, is
$$
X^k_t(i) (1 - X^k_t(j)) + X^k_t(j)(1 - X^k_t(i)).
$$
Therefore, on the complete graph, the probability that two different songs play at a meeting at time $t$ is
\begin{align*}
\sum_{i \neq j} {N \choose 2}^{-1} &\left( X^k_t(i) (1 - X^k_t(j)) + X^k_t(j)(1 - X^k_t(i)) \right) \\
&= \sum_{i \neq j} \frac{X^k_t(i) (1 - X^k_t(j))}{N (N - 1)} \\
&= \frac{N}{N - 1} M^k_t (1 - M^k_t) - \sum_i \frac{X^k_t(i)^2}{N(N - 1)} \\
&\geq M^k_t (1 - M^k_t) - \frac{1}{N - 1}.
\end{align*}

Recalling Lemma\nobreakspace \ref {M Steps}, at time $T_2$ we still have 
$$
M^k_{T_2} \in \left( \frac{1}{8 \sigma} - \frac{2 \eta}{N}, 1 - \frac{1}{8 \sigma} + \frac{2 \eta}{N} \right)
$$
and so at time $T_2$ we have
$$
M^k_{T_2} (1 - M^k_{T_2}) \geq \left( \frac{1}{8 \sigma} - \frac{2 \eta}{N}\right)^2
$$
Thus, the probability at time $T_2$ that fixation has occurred is bounded by
\begin{align*}
\Prb_{X(T_2)} ( T_{\fix} \geq 0 ) &\geq M^k_{T_2} (1 - M^k_{T_2}) - \frac{1}{N - 1} \\
&\geq \left( \frac{1}{8 \sigma} - \frac{2 \eta}{N}\right)^2 - \frac{1}{N - 1}.
\end{align*}
Applying the Strong Markov property, we can conclude that
\begin{align*}
\E T_{\fix} &\geq \E T_2 1(T_{\fix} \geq T_2) \\
&= \E T_2 \E ( 1(T_{\fix} \geq T_2) \vert T_2 ) \\
&= \E T_2 \left( \left( \frac{1}{8 \sigma} - \frac{2 \eta}{N}\right)^2 - \frac{1}{N - 1} \right) \\
&\geq \frac{N}{4  \eta^2 \sigma} \left( \left( \frac{1}{8 \sigma} - \frac{2 \eta}{N}\right)^2 - \frac{1}{N - 1} \right)  
\end{align*}
finishing the proof.

\end{proof}
\end{thm}

\section{Further Directions}

We conclude by presenting a few possible further directions for research on the iPod model.

\subsection{Improve the Fixation Time Bound}

Heuristically, from any initial configuration the processes $X^k_t$ mixes on a time scale of the order of the relaxation time $\lambda^{-1}$. Then, for any song $k$, when $x_t(i) \approx M_t$ we have $Q_t \approx M_t (1 - M_t)$ and so
$$
\E ( d M_t(1 - M_t) \vert F_{t-} ) \approx - \frac{2 \eta^2}{N} M_t( 1 - M_t) dt.
$$
Following through the same embedding and comparison arguments, we then find a fixation time of $O(N)$. Therefore we conjecture that for any initial configuration
$$
\E T_{\fix} = O(\lambda^{-1} + N) = O( \max(\lambda^{-1}, N ).
$$

\subsection{Remove the Dependence on $\sigma$}

When the processes $X^k_t$ are well mixed, i.e. when $x^k(i) \approx M^k$ again we have $Q_t \approx M_t (1 - M_t)$. Then, the $\sigma$-dimensional process $\{ M^k_t \}_{1 \leq k \leq \sigma, t \geq 0}$ has a comparable covariation structure to the $\sigma$-allele Wright Fisher Diffusion.

By a well known calculation \cite{durrett2002probability} the $\sigma$-allele process has an expected absorption time of $O(1)$, i.e. independent of $\sigma$. Therefore, we conjecture that by a similar embedding and comparison argument, the iPod process fixates in a time scale independent of the number of songs $\sigma$.

Combining this with our other conjectured improvements to \protect \MakeUppercase {T}heorem\nobreakspace \ref {Main Theorem}, we conclude with the following conjectured bound for the fixation time of the iPod model.

\begin{conjecture}{There exists a constant $C$ so that for any graph $\mathfrak{G}$ on $N$ vertices, the fixation time $T_{\fix}$ of the iPod model on $\mathfrak{G}$ with $\sigma$ songs, started from any initial configuration, satisfies
$$
\E T_{\fix} \leq \frac{C}{\eta (1 - \eta)} \max \left( N, \lambda^{-1} \right),
$$
where $\lambda$ is the spectral gap of $\mathfrak{G}$.
}
\end{conjecture}

%
%
%

\bibliographystyle{plain}

\bibliography{research}

\end{document}